# Estimation of the van de Vusse reactor via Carleman embedding

Dhruvi S. Bhatt* and Shambhu N. Sharma†

Electrical Engineering Department,
Sardar Vallabhbhai National Institute of Technology, Surat, Gujarat, India 395007
*bhattdhruvi427@gmail.com
†snsvolterra@gmail.com

**Abstract**

The van de Vusse reactor is an appealing benchmark problem in industrial control, since it has a non-minimum phase response. The van de Vusse stochasticity is attributed to the fluctuating input flow rate. The novelties of the paper are two. First, we utilize the surprising power of 'Itô's stochastic calculus for applications' to account for the van de Vusse stochasticity. Secondly, the Carleman embedding is unified with the Fokker-Planck equation for finding the estimation of the van de Vusse reactor. The revelation of the paper is that the Carleman linearized estimate of the van de Vusse reactor is more refined in contrast to the EKF-predicted estimate. This paper will be useful to practitioners aspiring for formal methods for stochastically perturbed nonlinear reactors as well as system theorists aspiring for applications of their theoretical results to practical problems.

**1. Introduction**

The Continuous Stirred-Tank Reactor (CSTR) is a significant element of major industrial chemical processes. In industrial control scenarios, a CSTR effectuating van de Vusse reaction is considered to be a benchmark control problem (van de Vusse, 1964; Engell & Klatt, 1993). Control studies are available in the deterministic framework of the van de Vusse reactor (Stack & Doyle, 1997; Perez et al., 2002). Moreover, Doyle et al. (1995) explored the Carleman linearization to arrive at the Volterra model for the non-linear van de Vusse reactor in order to design a superior control strategy. Despite the ubiquity of the van de Vusse reactor in the deterministic case, the stochastic case of the van de Vusse reactor is less researched. Hence, we restrict our discussion to the stochastic framework of the van de Vusse reactor.

Industrial chemical processes are often stochastically excited, subject to random fluctuations in input variables. The values of input flow rate, temperature and the reactant concentration undergo random fluctuations because of the change in ambient conditions, valve chatter attributed to interconnected sub-systems etc. (Aris & Amundson, 1958; Ligon & Amundson, 1981). Ignoring the randomness of these input variables will have an influence on 'the system parameter design and accurate capturing of qualitative characteristics of chemical reactors' (Berryman & Himmelblau, 1973). Further, the ignorance of the randomness will also lead to inaccurate estimates and stability issues. Moreover, the lack of prior assessments of the randomness in input variables on the system performance will lead to poor qualitative control of the process. Hence, it is evident to evaluate the effect of random fluctuations in input variables on the process performance from the state and parameter estimation viewpoints. A good source of state and parameter estimation perspectives of chemical reactors can be found in Dochain (2003). In Dochain (2003), non-linear observers were the subject of discussions. The observers hold for the deterministic case and the stochastic versions of the observers become predictors and filters. The notion of predictions is useful for the cases, where observations are not available. The notion of filters is useful, where observations are available. This necessitates the consideration of random fluctuations in input variables resulting in the stochastic non-linear framework of chemical reactors. That merits investigation for finding the estimation perspective. Despite the wide popularity of the van de Vusse reactor in industrial control, there is no formal, systematic and unified estimation theory available yet in literature for the 'van de Vusse reactor Itô stochastic differential equation'. The motivation for this paper comes from the fact that the older question in a new perspective unfolds unexplored insights (Holbrook & Bhatia, 2006).

The objective of this paper is to model and investigate the behaviour of Continuous Stirred Tank Reactor (CSTR) effectuating a van de Vusse reaction under fluctuations of the inlet flow rate. We achieve an analysis of a noise-influenced van de Vusse reactor. That involves two steps: Carleman embedding in the Itô framework and conditional moment evolutions via the Fokker-Planck equation. Khodadadi and Jazayeri-Rad (2011) achieved the noise analysis of the van de Vusse reactor in the Extended Kalman Filter (EKF) setting, where observations are available at discrete-time instants. In contrast to previously published results, the major ingredients of this paper are three celebrated results of applied mathematics, i.e. Carleman linearization, Fokker-Planck equation and the Itô theory. The idea of this paper resembles the following fact: recasting the finite-dimensional non-linear system into a system of infinite-dimensional linear systems (Carleman, 1932; Bellman & Richardson, 1963). That can be achieved via the Carleman embedding. The Fokker-Planck equation (Jazwinski, 1970; Sharma, 2008) is useful for the noise analysis of the SDEs, where no observations are available. The Itô framework allows a formal stochastic interpretation (Kunita, 2010). Thus, our approach unifies the Carleman linearization, the Fokker-Planck equation and the Itô calculus together for the van de Vusse reactor SDE. This paper exploits a different setting with the usefulness of new results. The Carleman setting of the paper offers a simpler realization of stochastic systems. That is attributed to fact that the Carleman linearization results bilinearization. Here, the 'Carleman embedding' means embedding the Carleman linearization into non-linear stochastic differential equations.

*A notational remark*: For a notational convenience, the time '$t$' is stated as an input argument of submatrices and a subscript of matrices respectively throughout the paper.

## 2. Carleman linearization of a stochastic van de Vusse reactor

A model of the van de Vusse reactor accounting for random fluctuations in the inlet flow rate can be described as (van de Vusse, 1964; Doyle et al., 1995)

$$\dot{C}_A = -k_1 C_A - k_3 C_A^2 + \frac{F_r}{v}(C_{Af} - C_A), \quad \dot{C}_B = k_1 C_A - k_2 C_B - \frac{F_r}{v} C_B. \tag{1}$$

Since the stochastic process is a collection of random variables as the time evolves, we consider the flow rate as a time-varying random variable. Thus, the stochastic evolution of the flow rate obeys the following SDE:

$$dF_r = -\alpha F_r dt + \beta dB_t. \tag{2}$$

For the simplified analysis, we consider the associated system parameters are constant coefficients of (1)-(2) evolution equations. The white noise has found its applications to model random forcing terms in dynamical systems and control. The white noise process has the zero correlation time as well as the white noise has an informal stochastic interpretation. As a result of this, the white noise is regarded as a generalized stochastic process in the theory of noise and stochastic processes (Wax, 1954).

Here, we explain why the Ornstein-Uhlenbeck (OU) process is accurate to model the flow rate of the reactor. In the state space setting, the OU process can be regarded as the output of an LTI system driven by the input Brownian motion. The current through the series RL circuit driven by the Brownian motion obeys the evolution of the OU process. The flow rate is analogous to the current. Thus, it is imperative to choose the flow rate as the OU stochastic process. The OU process obeys the real noise statistics. The real noise statistics has finite non-zero correlation time and the frequency-dependent spectrum. For these reasons, we consider the OU process statistics in lieu of the white noise statistics and the coloured noise case of the van de Vusse reactor is the subject of investigations. Rearranging the van de Vusse reactor in the SDE set up, we have

$$dx_t = f(t, x_t)dt + G(t, x_t)dB_t, \tag{3}$$

where

$$x_t = (x_1(t)\ x_2(t)\ x_3(t))^T = (C_A(t)\ C_B(t)\ F_r(t)),\ G(t, x_t) = (0\ 0\ \beta)^T,$$

$$f(t, x_t) = (f_1(t, x_t)\ f_2(t, x_t)\ f_3(t, x_t))^T$$

$$= (-k_1 C_A - k_3 C_A^2 + \frac{F_r}{v}(C_{Af} - C_A)\ \ k_1 C_A - k_2 C_B - \frac{F_r}{v} C_B\ \ -\alpha F_r)^T.$$

For the brevity of notations, the input argument $t$ associated with the components of the state vector $x_t$ is dropped. The Carleman realization to ordinary differential equations has received considerable attentions (Kawalski & Steeb 1991, p.75). However, the Carleman linearization to stochastic differential equations is relatively very scarce. Importantly, the Carleman linearization to non-linear stochastic differential equations leads to bilinear stochastic differential equations. In this paper, we do the Carleman linearization to the van de Vusse reactor SDE and we arrive at the van de Vusse bilinear SDE. Since the van de Vusse reactor (4) accounts for the square non-linearity, we consider the Carleman linearization order two. As a result of this, the augmented state vector accounting for the Carleman linearized state vector becomes $\begin{pmatrix} x_t \\ x_t^{(2)} \end{pmatrix}$, where

$$x_t^{(2)} = x_t \otimes x_t = (x_1 \ x_2 \ x_3)^T \otimes (x_1 \ x_2 \ x_3)^T,$$

and the notation $\otimes$ denotes the Kronecker product. The term $\binom{n+r-1}{r}$ denotes that the state $x_t^{(r)}$ has $\binom{n+r-1}{r}$ independent state variables, where $n$ and $r$ are the dimension of the state vector $x_t$ and the Kronecker power associated with the state respectively. Consider the case $r = 2$, then the state $x_t \in R^n$ and $x_t^{(2)} \in R^{\binom{n+1}{2}}$. In the specific case, we explain the dimension of the augmented state vector associated with the van de Vusse reactor SDE of the paper. For the SDE of the paper, consider $n = 3$, $r = 2$, we have six independent state variables associated with the state vector $x_t^{(2)}$. For the Carleman linearization order two, the dimension of the augmented state vector is nine. After accounting for independent state variables, we have

$$x_t^{(2)} = (x_1^2 \ x_1 x_2 \ x_1 x_3 \ x_2^2 \ x_2 x_3 \ x_3^2)^T, \ d\begin{pmatrix} x_t \\ x_t^{(2)} \end{pmatrix} = (dx_1 \ dx_2 \ dx_3 \ dx_1^2 \ dx_1 x_2 \ dx_1 x_3 \ dx_2^2 \ dx_2 x_3 \ dx_3^2)^T. \quad (4)$$

Using (3) and applying the stochastic differential rule for (4) (Karatzas & Shreve 1988, p. 154; Pugachev & Sinitsyn 1977, p. 163), we have the following evolution equations:

$$dx_1 = (-k_1 x_1 - k_3 x_1^2 + \frac{x_3}{v}(C_{Af} - x_1))dt, \ dx_2 = (k_1 x_1 - k_2 x_2 - \frac{x_3}{v} x_2)dt, \ dx_3 = (-\alpha x_3)dt + \beta dB_t, \quad (5a)$$

$$dx_1^2 = (-2k_1 x_1^2 - 2k_3 x_1^3 + \frac{2C_{Af}}{v} x_1 x_3 - \frac{2}{v} x_1^2 x_3)dt, \quad (5b)$$

$$dx_1 x_2 = (-(k_1 + k_2)x_1 x_2 + k_1 x_1^2 + \frac{C_{Af}}{v} x_2 x_3 - k_3 k x_1^2 x_2 - \frac{2}{v} x_1 x_2 x_3)dt, \quad (5c)$$

$$dx_1x_3 = (-(\alpha+k_1)x_1x_3 + \frac{C_{Af}}{v}x_3^2 - k_3x_1^2x_3 - \frac{1}{v}x_1x_3^2)dt + \beta x_1 dB_t, \quad (5d)$$

$$dx_2^2 = (2k_1x_1x_2 - 2k_2x_2^2 - \frac{2}{v}x_2^2x_3)dt, \quad dx_2x_3 = (-(\alpha+k_2)x_2x_3 + k_1x_1x_3 - \frac{1}{v}x_2x_3^2)dt + \beta x_2 dB_t, \quad (5e)$$

$$dx_3^2 = (-2\alpha x_3^2 + \beta^2)dt + 2\beta x_3 dB_t. \quad (5f)$$

Since the van de Vusse SDE accounts for the square non-linearity, we have the following: (i) the Carleman linearization order two (ii) the penultimate Kronecker power associated with the state vector is two. Here, we drop the explicit terms of the above system (5a)-(5f) of evolutions contributing to the higher-order than two, see Rugh (1970, p. 108). As a result, we arrive at the following system of stochastic differential equations:

$$dx_1 = (-k_1x_1 - k_3x_1^2 + \frac{x_3}{v}(C_{Af}-x_1))dt, \quad dx_2 = (k_1x_1 - k_2x_2 - \frac{x_3}{v}x_2)dt, \quad dx_3 = (-\alpha x_3)dt + \beta dB_t, \quad (6a)$$

$$dx_1^2 = (-2k_1x_1^2 + \frac{2C_{Af}}{v}x_3x_1)dt, \quad dx_1x_2 = (-(k_1+k_2)x_1x_2 + k_1x_1^2 + \frac{C_{Af}}{v}x_2x_3)dt, \quad (6b)$$

$$dx_1x_3 = (-(\alpha+k_1)x_1x_3 + \frac{C_{Af}}{v}x_3^2 + \beta x_1 dB_t)dt, \quad (6c)$$

$$dx_2^2 = (2k_1x_1x_2 - 2k_2x_2^2)dt, \quad dx_2x_3 = (-(\alpha+k_2)x_2x_3 + k_1x_1x_3 + \beta x_2 dB_t)dt, \quad (6d)$$

$$dx_3^2 = (-2\alpha x_3^2 + \beta^2)dt + 2\beta x_3 dB_t. \quad (6e)$$

It is important to mention that the order of the non-linearity decides the Carleman linearization order of non-linear functions. The Carleman linearization order contributes to the dimension of the augmented state vector. Suppose the Carleman linearization order is $N$, then the dimension of the augmented state vector is $\sum_{1 \leq k \leq N} n^k$. After accounting for independent state variables in the state vector, the dimension of the augmented state vector is $\sum_{1 \leq r \leq N} \binom{n+r-1}{r}$. Rearranging the set (6a)-(6e) in the SDE set up, we are led to the following bilinear SDE:

$$d\xi_t = (A_0(t) + A_t\xi_t)dt + D_t\xi_t dB_t + G_t dB_t, \quad (7)$$

where the state vector $\xi_t = (\xi_i) = (x_t^{(k)})$, and $k$ is the Kronecker power, $1 \leq k \leq N$ $1 \leq i \leq \sum_{1 \leq k \leq N}\binom{n+k-1}{k}$. Note that the augmented state vector $\xi_t \in R^{\sum_{1 \leq k \leq N}\binom{n+k-1}{k}}$. In the partitioned matrix-vector format, the van de Vusse Carleman linearized reactor SDE becomes

$$d\begin{pmatrix}x_t\\x_t^{(2)}\end{pmatrix}=\left(\begin{pmatrix}A_{01}(t)\\A_{02}(t)\end{pmatrix}+\begin{pmatrix}A_{11}(t)&A_{12}(t)\\A_{21}(t)&A_{22}(t)\end{pmatrix}\begin{pmatrix}x_t\\x_t^{(2)}\end{pmatrix}\right)dt+\left(\begin{pmatrix}D_{11}(t)&D_{12}(t)\\D_{21}(t)&D_{22}(t)\end{pmatrix}\begin{pmatrix}x_t\\x_t^{(2)}\end{pmatrix}+\begin{pmatrix}G_1(t)\\G_2(t)\end{pmatrix}\right)dB_t. \quad (8)$$

For the time-invariant case, we have

$$A_{11}(t)=\begin{pmatrix}-k_1 & 0 & \frac{C_{Af}}{v}\\ k_1 & -k_2 & 0\\ 0 & 0 & -\alpha\end{pmatrix},\ A_{12}(t)=\begin{pmatrix}-k_3 & 0 & -\frac{1}{v} & 0 & 0 & 0\\ 0 & 0 & 0 & 0 & -\frac{1}{v} & 0\\ 0 & 0 & 0 & 0 & 0 & 0\end{pmatrix}, \quad (9a)$$

$$A_{22}(t)=\begin{pmatrix}-2k_1 & 0 & \frac{2C_{Af}}{v} & 0 & 0 & 0\\ k_1 & -(k_2+k_2) & 0 & 0 & \frac{C_{Af}}{v} & 0\\ 0 & 0 & -(\alpha+k_1) & 0 & 0 & \frac{C_{Af}}{v}\\ 0 & 2k_1 & 0 & -2k_2 & 0 & 0\\ 0 & 0 & k_1 & 0 & -(\alpha+k_2) & 0\\ 0 & 0 & 0 & 0 & 0 & -2\alpha\end{pmatrix},\ D_{21}(t)=\begin{pmatrix}0 & 0 & 0\\ 0 & 0 & 0\\ \beta & 0 & 0\\ 0 & 0 & 0\\ 0 & \beta & 0\\ 0 & 0 & 2\beta\end{pmatrix}, \quad (9b)$$

$$G_1(t)=(0\ 0\ \beta)^T, G_2(t)=(0\ 0\ 0\ 0\ 0\ 0)^T, A_{01}(t)=(0\ 0\ 0)^T, A_{02}(t)=(0\ 0\ 0\ 0\ 0\ \beta^2)^T. \quad (9c)$$

Furthermore, we adopt the notational brevity for zero matrices, i.e. $A_{21}(t)=(A_{21}^{ij})$, where $1\leq i\leq 6, 1\leq j\leq 3$ and $A_{21}^{ij}(t)=0$. Note that $D_{11}(t)=(D_{11}^{ij})$, where $1\leq i\leq 3, 1\leq j\leq 3$ and $D_{11}^{ij}(t)=0$, $D_{12}(t)=(D_{12}^{ij})$, where $1\leq i\leq 3, 1\leq j\leq 6$ and $D_{12}^{ij}(t)=0$, $D_{22}(t)=(D_{22}^{ij})$, where $1\leq i\leq 6, 1\leq j\leq 6$ and $D_{22}^{ij}(t)=0$.

In the general setting, the dimension of the state vector $x_t$ is $n$ and the dimension of the state vector $x_t^{(2)}$ is $\frac{n(n+1)}{2}$. Thus, the augmented state vector has the size $n+\frac{n(n+1)}{2}$. The sizes of the principal submatrices $A_{11}(t)$ and $A_{22}(t)$ are $n\times n$ and $\frac{n(n+1)}{2}\times\frac{n(n+1)}{2}$. The sizes of the submatrices $A_{12}(t)$ and $A_{21}(t)$ are $n\times\frac{n(n+1)}{2}$ and $\frac{n(n+1)}{2}\times n$. The conditional mean and variance evolutions of (A.13)-(A.14) for the Carleman linearized Itô bilinear SDE can be recast as

$$d\widehat{x}_t = (A_{01}(t) + A_{11}(t)\widehat{x}_t + A_{12}(t)\widehat{x}_t^{(2)})dt, \quad d\widehat{x}_t^{(2)} = (A_{02}(t) + A_{21}(t)\widehat{x}_t + A_{22}(t)\widehat{x}_t^{(2)})dt, \tag{10}$$

$$\frac{dP_{x_t x_t}}{dt} = (P_{x_t x_t} A_{11}^T(t) + P_{x_t x_t^{(2)}} A_{12}^T(t) + A_{11}(t)P_{x_t x_t} + A_{12}(t)P_{x_t^{(2)} x_t} + G_1(t)G_1^T(t) + D_{11}(t)\widehat{x}_t G_1^T(t) + D_{12}(t)\widehat{x}_t^{(2)} G_1^T(t) + G_1(t)\widehat{x}_t^T D_{11}^T(t)$$

$$+ G_1(t)\widehat{x}_t^{(2)^T} D_{12}^T(t) + D_{11}(t)(P_{x_t x_t} D_{11}^T(t) + P_{x_t x_t^{(2)}} D_{12}^T(t)) + D_{12}(t)(P_{x_t^{(2)} x_t} D_{11}^T(t) + P_{x_t^{(2)} x_t^{(2)}} D_{12}^T(t))$$

$$+ D_{11}(t)(\widehat{x}_t \widehat{x}_t^T D_{11}^T(t) + \widehat{x}_t \widehat{x}_t^{(2)^T} D_{12}^T(t)) + D_{12}(t)(\widehat{x}_t^{(2)} \widehat{x}_t^T D_{11}^T(t) + \widehat{x}_t^{(2)} \widehat{x}_t^{(2)^T} D_{12}^T(t))), \tag{11a}$$

$$\frac{dP_{x_t x_t^{(2)}}}{dt} = (P_{x_t x_t} A_{21}^T(t) + P_{x_t x_t^{(2)}} A_{22}^T(t) + A_{11}(t)P_{x_t x_t^{(2)}} + A_{12}(t)P_{x_t^{(2)} x_t^{(2)}} + G_1(t)G_2^T(t) + D_{11}(t)\widehat{x}_t G_2^T(t) + D_{12}(t)\widehat{x}_t^{(2)} G_2^T(t)$$

$$+ G_1(t)\widehat{x}_t^T D_{21}^T(t) + G_1(t)\widehat{x}_t^{(2)^T} D_{22}^T(t) + D_{11}(t)(P_{x_t x_t} D_{21}^T(t) + P_{x_t x_t^{(2)}} D_{22}^T(t))$$

$$+ D_{12}(t)(P_{x_t^{(2)} x_t} D_{21}^T(t) + P_{x_t^{(2)} x_t^{(2)}} D_{22}^T(t)) + D_{11}(t)(\widehat{x}_t \widehat{x}_t^T D_{21}^T(t) + \widehat{x}_t \widehat{x}_t^{(2)^T} D_{22}^T(t))$$

$$+ D_{12}(t)(\widehat{x}_t^{(2)} \widehat{x}_t^T D_{21}^T(t) + \widehat{x}_t^{(2)} \widehat{x}_t^{(2)^T} D_{22}^T(t))), \tag{11b}$$

$$\frac{dP_{x_t^{(2)} x_t}}{dt} = (P_{x_t^{(2)} x_t} A_{11}^T(t) + P_{x_t^{(2)} x_t^{(2)}} A_{12}^T(t) + A_{21}(t)P_{x_t x_t} + A_{22}(t)P_{x_t^{(2)} x_t} + G_2(t)G_1^T(t) + D_{21}(t)\widehat{x}_t G_1^T(t) + D_{22}(t)\widehat{x}_t^{(2)} G_1^T(t)$$

$$+ G_2(t)\widehat{x}_t^T D_{11}^T(t) + G_2(t)\widehat{x}_t^{(2)^T} D_{12}^T(t) + D_{21}(t)(P_{x_t x_t} D_{11}^T(t) + P_{x_t x_t^{(2)}} D_{12}^T(t))$$

$$+ D_{22}(t)(P_{x_t^{(2)} x_t} D_{11}^T(t) + P_{x_t^{(2)} x_t^{(2)}} D_{12}^T(t)) + D_{21}(t)(\widehat{x}_t \widehat{x}_t^T D_{11}^T(t) + \widehat{x}_t \widehat{x}_t^{(2)^T} D_{12}^T(t))$$

$$+ D_{22}(t)(\widehat{x}_t^{(2)} \widehat{x}_t^T D_{11}^T(t) + \widehat{x}_t^{(2)} \widehat{x}_t^{(2)^T} D_{12}^T(t))), \tag{11c}$$

$$\frac{dP_{x_t^{(2)} x_t^{(2)}}}{dt} = (P_{x_t^{(2)} x_t} A_{21}^T(t) + P_{x_t^{(2)} x_t^{(2)}} A_{22}^T(t) + A_{21}(t)P_{x_t x_t^{(2)}} + A_{22}(t)P_{x_t^{(2)} x_t^{(2)}} + G_2(t)G_2^T(t) + D_{21}(t)\widehat{x}_t G_2^T(t) + D_{22}(t)\widehat{x}_t^{(2)} G_2^T(t)$$

$$+ G_2(t)\widehat{x}_t^T D_{21}^T(t) + G_1(t)\widehat{x}_t^{(2)^T} D_{22}^T(t) + D_{21}(t)(P_{x_t x_t} D_{21}^T(t) + P_{x_t x_t^{(2)}} D_{22}^T(t))$$

$$+ D_{22}(t)(P_{x_t^{(2)} x_t} D_{21}^T(t) + P_{x_t^{(2)} x_t^{(2)}} D_{22}^T(t)) + D_{21}(t)(\widehat{x}_t \widehat{x}_t^T D_{21}^T(t) + \widehat{x}_t \widehat{x}_t^{(2)^T} D_{22}^T(t))$$

$$+ D_{22}(t)(\widehat{x}_t^{(2)} \widehat{x}_t^T D_{21}^T(t) + \widehat{x}_t^{(2)} \widehat{x}_t^{(2)^T} D_{22}^T(t))). \tag{11d}$$

It is worth to mention that the above coupled conditional mean and conditional variance evolution equations hold for the Itô bilinear stochastic differential equation for the Carleman linearization order two and the arbitrary dimension of the state vector $x_t$. The intermediate steps associated with (11a)-(11d) can be found in the *appendix*, see (A.14). It is important to note that the notation $\widehat{x}_t^{(2)}$ denotes the action of the conditional expectation operator on the term $x_t^{(2)}$ in the Kronecker product setting (Bellman, 1960; Holbrook & Bhatia, 2006). The term $\widehat{x}_t^{(2)^T}$ denotes the transpose of the term $\widehat{x}_t^{(2)}$.

**Remark 1.** The Carleman linearization order decides the partitioning of the conditional mean vector and the conditional variance matrix. The conditional mean (10) and the conditional variance evolution (11a)-(11d) are associated with (8). State submatrices and the process noise coefficient submatrices of the van de Vusse reactor SDE are given in (9a)-(9c). The proof of conditional mean and variance evolutions can be traced back to the Fokker-Planck equation. The van de Vusse non-linear SDE can be rephrased as a bilinear Itô stochastic differential equation with the augmented state vector. That is a consequence of the Carleman linearization coupled with the celebrated Itô stochastic differential.

**Remark 2.** Since we wish to analyze the van de Vusse SDE via the Carleman embedding, the Fokker-Planck equation of the van de Vusse 'bilinear' SDE is the cornerstone. The intermediate steps 'from the Fokker-Planck equation to the conditional moment evolutions' for the bilinear Itô SDE with the augmented state can be found in the *appendix*. To write the Fokker-Planck equation, we rephrase the augmented state vector $\xi_t$ of the van de Vusse SDE, i.e. $\xi_t = (\xi_i) = (x_t^{(k)})$, where $1 \leq k \leq 2$, $1 \leq i \leq 6$.

Alternatively $\xi_t = (x_1 \ x_2 \ x_3 \ x_1^2 \ x_1 x_2 \ x_1 x_3 \ x_2^2 \ x_2 x_3 \ x_3^2)^T$. Note that the entries of the augmented state vector $\xi_t$ are the scalar and obey the property of phase variables. Thus, the stochastic van de Vusse reactor with square non-linearity has the augmented state vector $\xi_t \in R^{n + \frac{n(n+1)}{2}}$, $n = 3$. For the brevity of notations, we choose the conditional probability density notation

$$p(\xi,t) = p(x_1, x_2, x_3, x_1^2, x_1 x_2, x_1 x_3, x_2^2, x_2 x_3, x_3^2, t). \quad (12)$$

Thus

$$\frac{\partial p(\xi,t)}{\partial t} = -\frac{\partial}{\partial x_1}(-k_1 x_1 + \frac{C_{Af}}{v} x_3 - k_3 x_1^2 - \frac{1}{v} x_1 x_3) p - \frac{\partial}{\partial x_2}(k_1 x_1 - k_2 x_2 - \frac{1}{v} x_2 x_3) p + \frac{\partial}{\partial x_3} \alpha x_3 p$$

$$-\frac{\partial}{\partial x_1^2}(-2 k_1 x_1^2 + \frac{2 C_{Af}}{v} x_1 x_3) p - \frac{\partial}{\partial x_1 x_2}(k_1 x_1^2 - (k_1 + k_2) x_1 x_2 + \frac{C_{Af}}{v} x_2 x_3) p$$

$$-\frac{\partial}{\partial x_1 x_3}(-(\alpha + k_1) x_1 x_3 + \frac{C_{Af}}{v} x_3^2) p - \frac{\partial}{\partial x_2^2}(2 k_1 x_1 x_2 - 2 k_1 x_2^2) p$$

$$-\frac{\partial}{\partial x_2 x_3}(k_1 x_1 x_3 - (\alpha + k_2) x_2 x_3) p + \frac{\partial}{\partial x_3^2} 2 \alpha x_3^2 p$$

$$+\frac{1}{2} (\frac{\partial^2}{\partial x_3 \partial x_3} \beta^3 p + \frac{\partial^2}{\partial x_3 \partial x_1 x_3} \beta^2 x_1 p + \frac{\partial^2}{\partial x_3 \partial x_2 x_3} \beta^2 x_2 p + \frac{\partial^2}{\partial x_3 \partial x_3^2} 2 \beta^2 x_3 p$$

$$+\frac{\partial^2}{\partial x_1 x_3 \partial x_3} \beta^2 x_1 p + \frac{\partial^2}{\partial x_1 x_3 \partial x_1 x_3} \beta^2 x_1^2 p + \frac{\partial^2}{\partial x_1 x_3 \partial x_2 x_3} \beta^2 x_1 x_2 p + \frac{\partial^2}{\partial x_1 x_3 \partial x_3^2} 2 \beta^2 x_1 x_3 p$$

$$+\frac{\partial^2}{\partial x_2 x_3 \partial x_3}\beta^2 x_2 p + \frac{\partial^2}{\partial x_2 x_3 \partial x_1 x_3}\beta^2 x_1 x_2 p + \frac{\partial^2}{\partial x_2 x_3 \partial x_2 x_3}\beta^2 x_2^2 p + \frac{\partial^2}{\partial x_2 x_3 \partial x_3^2}2\beta^2 x_2 x_3 p$$

$$+\frac{\partial^2}{\partial x_3^2 \partial x_3}2\beta^2 x_3 p + \frac{\partial^2}{\partial x_3^2 \partial x_1 x_3}2\beta^2 x_1 x_3 p + \frac{\partial^2}{\partial x_3^2 \partial x_2 x_3}2\beta^2 x_2 x_3 p + \frac{\partial^2}{\partial x_3^2 \partial x_3^2}4\beta^2 x_3^2 p). \tag{13}$$

### 3. Conditional moment evolutions of the van de Vusse reactor SDE

Here, we write the conditional moment evolutions of the Carleman linearized van de Vusse SDE, see (10). More specifically, the element-wise conditional mean evolution equation of the bilinear van de Vusse SDE, which is a consequence of (10)-(11d), boils down to

$$d\hat{x}_1 = (-k_1\hat{x}_1 + \frac{C_{Af}}{v}\hat{x}_3 - k_3 P_{x_1} - k_3\hat{x}_1^2 - \frac{1}{v}P_{x_1 x_3} - \frac{\hat{x}_1 \hat{x}_3}{v})dt, \tag{14a}$$

$$d\hat{x}_2 = (k_1\hat{x}_1 - k_2\hat{x}_2 - \frac{1}{v}P_{x_2 x_3} - \frac{\hat{x}_2 \hat{x}_3}{v})dt, \quad d\hat{x}_3 = -\alpha\hat{x}_3 dt, \tag{14b}$$

$$d\hat{x}_1^2 = (-2k_1 P_{x_1} - 2k_1\hat{x}_1^2 + \frac{2C_{Af}}{v}P_{x_1 x_3} + \frac{2C_{Af}}{v}\hat{x}_1\hat{x}_3)dt, \tag{14c}$$

$$d\widehat{x_1 x_2} = (k_1 P_{x_1} + k_1\hat{x}_1^2 - (k_1+k_2)P_{x_1 x_2} - (k_1+k_2)\hat{x}_1\hat{x}_2 + \frac{C_{Af}}{v}P_{x_2 x_3} + \frac{C_{Af}}{v}\hat{x}_2\hat{x}_3)dt, \tag{14d}$$

$$d\widehat{x_1 x_3} = (-(\alpha+k_1)P_{x_1 x_3} - (\alpha+k_1)\hat{x}_1\hat{x}_3 + \frac{C_{Af}}{v}P_{x_3} + \frac{C_{Af}}{v}\hat{x}_3^2)dt, \tag{14e}$$

$$d\hat{x}_2^2 = (2k_1 P_{x_1 x_2} + 2k_1\hat{x}_1\hat{x}_2 - 2k_2 P_{x_2} - 2k_2\hat{x}_2^2)dt, \tag{14f}$$

$$d\widehat{x_2 x_3} = (k_1 P_{x_1 x_3} + k_1\hat{x}_1\hat{x}_3 - (\alpha+k_2)P_{x_2 x_3} - (\alpha+k_2)\hat{x}_2\hat{x}_3)dt, \quad d\hat{x}_3^2 = (\beta^2 - 2\alpha P_{x_3} - 2\alpha\hat{x}_3^2)dt. \tag{14g}$$

An elementary proof of (14a)-(14g) can be sketched using the notation of (12) and the van de Vusse Fokker-Planck equation of (13). It is natural to ask how the conditional variance terms evolve, which are associated with the conditional mean evolutions. The conditional variance evolutions for the van de Vusse SDE become a special case of the coupled evolution equations (A.11)-(A.12). Thus,

$$dP_{x_1} = d(\hat{x}_1^2 - \hat{x}_1^2) = d\hat{x}_1^2 - d\hat{x}_1^2$$

$$= (-2k_1 P_{x_1} + \frac{2C_{Af}}{v}P_{x_1 x_3} + 2k_3\hat{x}_1 P_{x_1} + 2k_3\hat{x}_1^3 + \frac{2}{v}\hat{x}_1 P_{x_1 x_3} + \frac{2\hat{x}_1^2 \hat{x}_3}{v})dt, \tag{15a}$$

$$dP_{x_2} = d(\hat{x}_2^2 - \widehat{x_2^2}) = d\widehat{x_2^2} - d\hat{x}_2^2 = (2k_1 P_{x_1 x_2} - 2k_2 P_{x_2} + \frac{2}{v}\hat{x}_2 P_{x_2 x_3} + \frac{2\hat{x}_2^2 \hat{x}_3}{v})dt, \qquad (15b)$$

$$dP_{x_3} = (\beta^2 - 2\alpha P_{x_3})dt, \qquad (15c)$$

$$dP_{x_1 x_2} = d(\widehat{x_1 x_2}) - \hat{x}_1 d\hat{x}_2 - \hat{x}_2 d\hat{x}_1 - d\hat{x}_1 d\hat{x}_2$$

$$= (k_1 P_{x_1} + k_3 \hat{x}_2 P_{x_1} - (k_1 + k_2) P_{x_1 x_2} + \frac{1}{v}\hat{x}_2 P_{x_1 x_3} + \frac{C_{Af}}{v} P_{x_2 x_3} + \frac{1}{v}\hat{x}_1 P_{x_2 x_3} + k_3 \hat{x}_1^2 \hat{x}_2 + \frac{2\hat{x}_1 \hat{x}_2 \hat{x}_3}{v})dt, \qquad (15d)$$

$$dP_{x_1 x_3} = d(\widehat{x_1 x_3}) - \hat{x}_1 d\hat{x}_3 - \hat{x}_3 d\hat{x}_1 - d\hat{x}_1 d\hat{x}_3$$

$$= (-(\alpha + k_1) P_{x_1 x_3} + \frac{C_{Af}}{v} P_{x_3} + k_3 \hat{x}_3 P_{x_1} + k_3 \hat{x}_1^2 \hat{x}_3 + \frac{1}{v}\hat{x}_3 P_{x_1 x_3} + \frac{\hat{x}_1 \hat{x}_3^2}{v})dt, \qquad (15e)$$

$$dP_{x_2 x_3} = d(\widehat{x_2 x_3}) - \hat{x}_2 d\hat{x}_3 - \hat{x}_3 d\hat{x}_2 - d\hat{x}_2 d\hat{x}_3 = (k_1 P_{x_1 x_3} - (\alpha + k_2) P_{x_2 x_3} + \frac{1}{v}\hat{x}_3 P_{x_2 x_3} + \frac{\hat{x}_2 \hat{x}_3^2}{v})dt. \qquad (15f)$$

## 4. Numerical Simulations

To test the effectiveness of the proposed estimation method, i.e. the Carleman linearization-based, we perform numerical simulations of the van de Vusse reactor. The simulation is carried out for two different sets of initial conditions and operating parameters associated with the van de Vusse reactor. The performance of the proposed estimation method for the set of values is evaluated by contrasting it with the benchmark estimation method, namely Extended Kalman Filter-based prediction. The Carleman linearization is carried out at the operating points mentioned in Table 1 (Doyle et al., 1995).

The proposed estimation method of this paper is implemented in MATLAB© on Intel(R) Core(TM) i5-5200U laptop CPU clocked at 2.20 GHz with 8.00GB RAM. It is difficult to choose analytically a good estimate of the Carleman linearization order (Bellman & Richardson, 1963), see Bellman and Richardson (1961) as well. Note that this paper considers the order of the Carleman linearization two and the dimension of the state vector is three. We demonstrate graphically the usefulness of the Carleman linearization order two. Consider the first set of initial conditions and system parameters, i.e. $x_1(0) = 3, x_2(0) = 1.12, x_3(0) = 0.009528, \alpha = 0.1, \beta = 0.044$ and $P_{x_3}(0) = 0.01$. Here, the terms $\alpha$ and $\beta$ are the damping factor and the noise co-efficient of the OU process respectively.

**Table 1**
The first set of parameters

| Parameters | Values | Units |
|---|---|---|
| $k_1$ | 0.01388 | $s^{-1}$ |
| $k_2$ | 0.02778 | $s^{-1}$ |
| $k_3$ | 0.002778 | $lmol^{-1}s^{-1}$ |
| $C_{Af}$ | 0.0027 | $moll^{-1}$ |
| $v$ | 10 | $L$ |
| $C_A$ | 3 | $moll^{-1}$ |
| $C_B$ | 1.12 | $moll^{-1}$ |
| $F_r$ | 0.009528 | $ls^{-1}$ |

Fig. 1 shows the comparison of Carleman linearized state trajectories with the true state trajectories, i.e. non-linear stochastic differential equations. Fig. 1 has two parts, Fig. 1(a) and Fig. 1(b). Fig. 1(a) displays the state $x_1$ trajectory corresponding to the concentration of the reactant $A$. The solid line of Fig. 1(a) illustrates the state $x_1$ evolution, which is a consequence of the exact van de Vusse SDE. The dash-dash line of Fig. 1(a) shows the state $x_1$ evolution, which is a consequence of the Carleman linearized SDE. Similar interpretations hold for Fig. 1(b) associated with the concentration of reactant $B$.

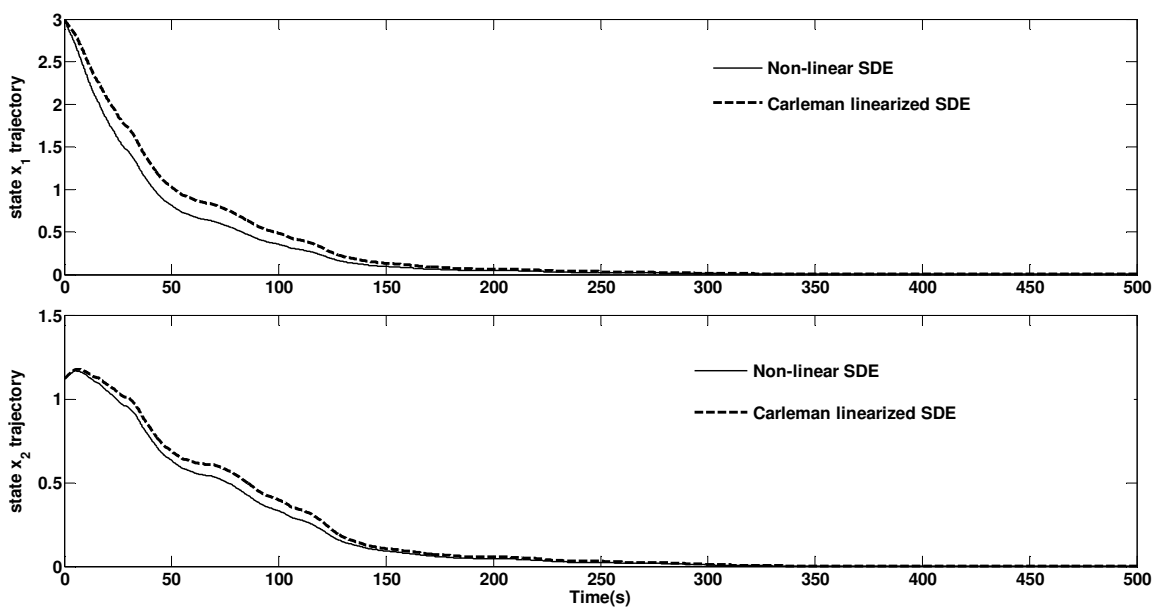

**Fig. 1.** A comparison between the van de Vusse SDE and its Carleman linearized SDE.

Fig. 1 reveals that the Carleman linearized trajectories associated with the states $x_1$ and $x_2$ track the non-linear state trajectories. The closeness of the linearized trajectories to that of the non-linear noisy trajectories unfolds that the Carleman linearization of the order two for the van de Vusse reactor captures the square nonlinearity quite well.

Fig. 2 has two parts, Fig. 2(a) and Fig. 2(b). Fig. 2(a) contrasts the Carleman linearized estimated state trajectory, the Extended Kalman Filter-predicted trajectory and the true state trajectory. In Figs. 2(a) and 2(b) solid lines indicate the true state trajectories, dash-dash and dot-dot lines denotes the Carleman linearized estimate and EKF-predicted state trajectories respectively. The three trajectories of Fig. 2(a) are associated with the state $x_1$. the prediction error evolution exploiting the Carleman linearized estimated state. The dash-dash line illustrates the EKF-predicted estimate error. Note that the EKF-predicted estimate is an estimate associated with 'between the observations' case of the EKF for the continuous state-discrete measurement. That is a non-linear predicted estimate unaccounting for observations (Jazwinski 1970, p. 278). Note that the Carleman linearized estimate of the paper is also a predicted estimate, since the observations are not available. In the absence of observations, the estimates are predicted (Jazwinski 1970, p. 179). The Carleman linearized estimated state trajectories are the consequence of the conditional mean equations of the van de Vusse reactor, see (14a)-(14g), in combination with the initial conditions and the operating parameters. Note that the true state trajectories are a consequence of the non-linear SDE model of the van de Vusse reactor, see (3)-(4) and the values of Table 1. The three trajectories of Fig. 2(b) are associated with the state $x_2$.

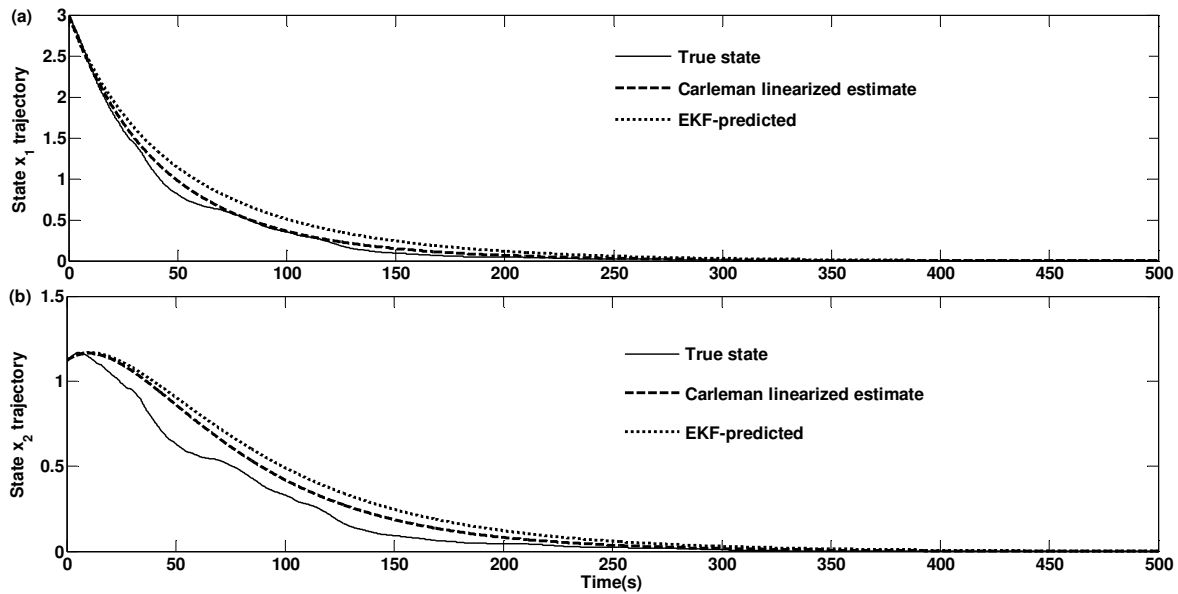

**Fig. 2.** A comparison between three trajectories.

Fig. 2 reveals that the Carleman linearized estimated trajectories for the states successfully track the non-linear stochastic state trajectories under the influence of random fluctuations in the inlet flow rate. Moreover, the Carleman linearized estimated state trajectories are much closer to the true state trajectories in contrast to the EKF-predicted trajectories.

Figs. 3(a) and 3(b) show a comparison between two absolute prediction error evolutions for the states $x_1$ and $x_2$ respectively. The solid line illustrates the prediction error evolution exploiting the Carleman linearized estimated state. The dash-dash line illustrates the EKF-predicted estimate error. Note that the absolute prediction error $e_t$ is defined as an absolute difference between the true state and the mean values of the state, i.e.

$$e_t = |x_t - \widehat{x}_t|, \; t = 1, 2.$$

It is quite evident from the absolute prediction error trajectory, displayed in Fig. 3, that the absolute prediction error associated with the Carleman linearized estimate has the less value in contrast to the EKF-prediction. For the given set of system parameters in Table 1, the maximum absolute prediction error associated with the Carleman linearized estimate of the state $x_1$ is less than $0.20$, on the other hand, the maximum absolute error associated with the EKF-predicted estimate is greater than $0.30$, see Fig. 3(a). Fig. 3(b) illustrates the maximum absolute prediction error associated with the Carleman linearized estimate is relatively less than $0.25$, on the other hand, the maximum absolute error with the EKF-predicted estimate is quite more than $0.25$. The less value of the absolute prediction error, associated with the proposed method, concords the greater closeness of the estimated trajectories via the Carleman linearization with true state trajectories, see Fig. 2 as well.

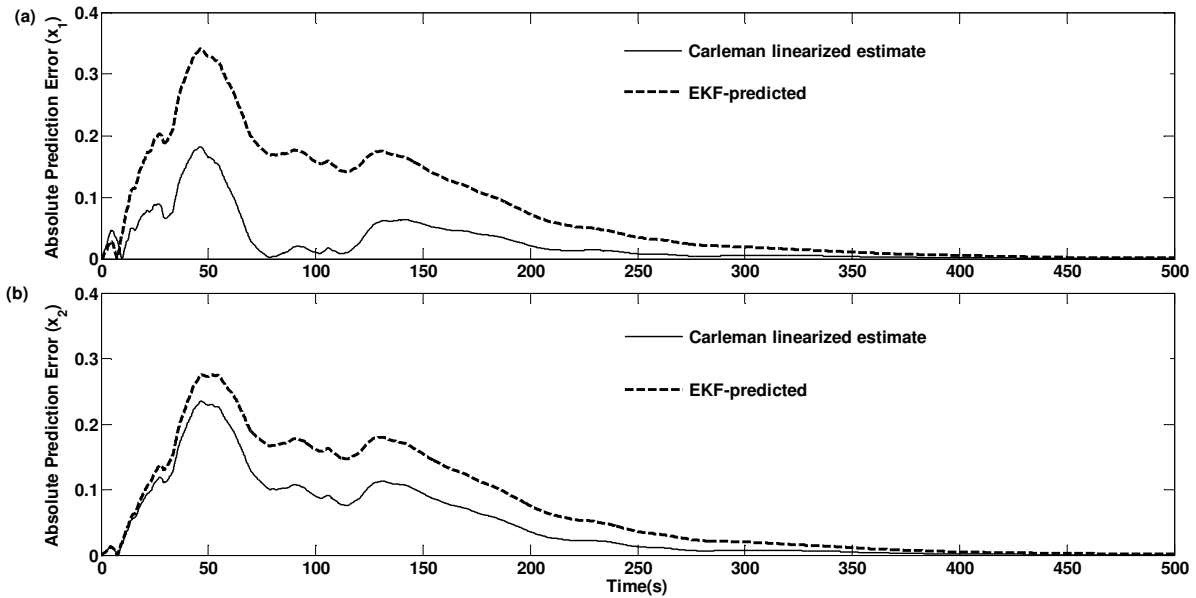

**Fig. 3.** Absolute error comparison for the both states.

The expected value of the square of the prediction error gives the variance trajectories, see Fig. 4. The variance trajectory tells random fluctuations in the mean trajectory. The conditional variance trajectories for the both states of the van de Vusse reactor are shown in Figs. 4(a) and 4(b) respectively. Note that the conditional variance trajectories associated with the Carleman estimation are a consequence of conditional variance evolution equations derived theoretically, see (15a)-(15f). In comparison to conditional variance trajectories of the EKF-prediction, the Carleman linearization-based conditional variance trajectories show quite less variance, see Fig. 4. Note that Fig. 4 has two parts, Fig. 4(a) and Fig. 4(b). Table 2 shows the values of conditional variance, associated with the first set of parameters, at several time instances. The less variance suggests the better estimate. That is indicative of less random fluctuations in the most probable trajectory.

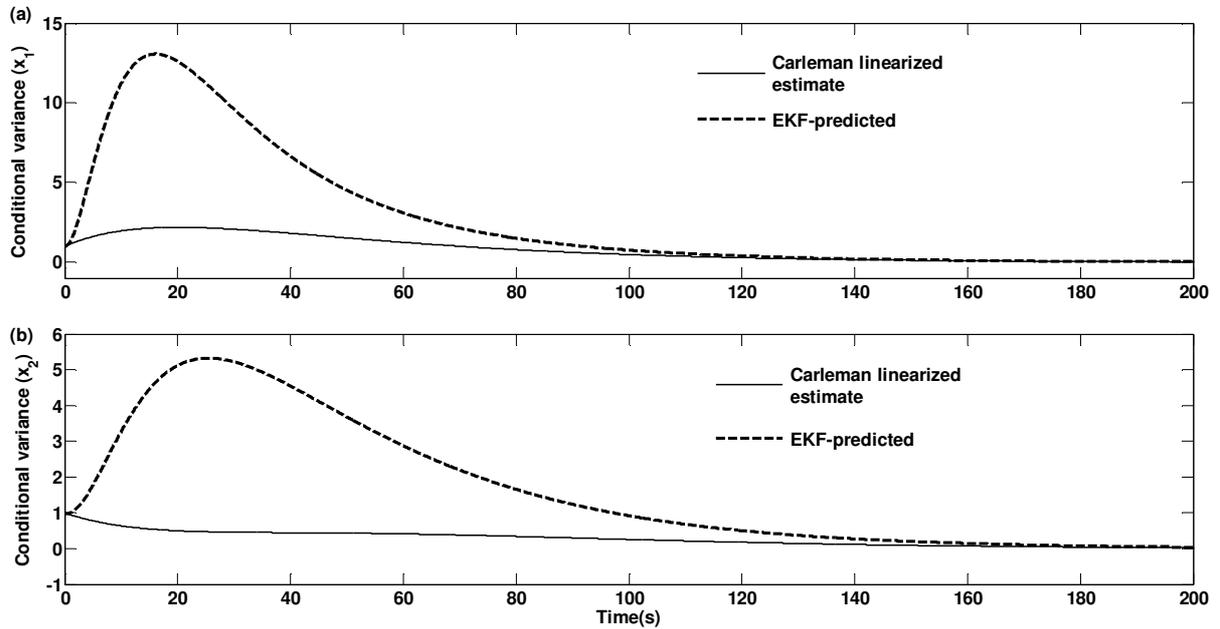

**Fig. 4.** Conditional variance trajectories of the states $x_1$ and $x_2$.

Table 2 unfolds the Carleman linearized estimate is better than that of the EKF-predicted estimate. Hence, the proposed method is a better performer to estimate the state trajectories of an isothermal van de Vusse reactor under OU process-driven random fluctuations considered in the inlet flow rate. It is quite evident from the numerical simulation that the proposed estimation method, i.e. the Carleman linearization based, is superior to the benchmark estimation method, i.e. the EKF-predicted.

**Table 2**
Conditional variance of concentrations for the first set of parameters

| Time | Carleman linearized estimated state $x_1$ | EKF-predicted state $x_1$ | Carleman linearized estimated state $x_2$ | EKF-predicted state $x_2$ |
|---|---|---|---|---|
| 0.5 | 1.08 | 1.09 | 0.97 | 0.99 |
| 5 | 1.62 | 6.10 | 0.77 | 1.78 |
| 10 | 1.97 | 11.23 | 0.63 | 3.28 |
| 20 | 2.19 | 12.53 | 0.50 | 5.13 |
| 50 | 1.52 | 4.53 | 0.43 | 3.68 |
| 100 | 0.48 | 0.76 | 0.25 | 0.92 |
| 150 | 1.13 | 0.16 | 0.09 | 0.20 |
| 200 | 0.03 | 0.04 | 0.03 | 0.05 |

Now consider the second set of operating parameters mentioned in Table 3 (Åkesson & Toivonen, 2006). The initial conditions associated with the operating parameters of Table 3 are $x_1(0) = 1.235, x_2(0) = 1, x_3(0) = 0.0152, \alpha = 0.01, \beta = 0.044$ and $P_{x_3}(0) = 0.09$.

**Table 3**
The second set of operating parameters

| Parameters | Values | Units |
|---|---|---|
| $k_1$ | 0.0141 | $s^{-1}$ |
| $k_2$ | 0.0141 | $s^{-1}$ |
| $k_3$ | 0.00187 | $lmol^{-1}s^{-1}$ |
| $C_{Af}$ | 0.00141 | $moll^{-1}$ |
| $v$ | 10 | $L$ |
| $C_A$ | 1.235 | $moll^{-1}$ |
| $C_B$ | 1 | $moll^{-1}$ |
| $F_r$ | 0.0152 | $ls^{-1}$ |

Numerical simulations utilizing the second set of parameters adopt a similar procedure to that of the first set of parameters. Fig. 5(a) show the comparison between the true state trajectory, Carleman linearized estimated state trajectory and the Extended Kalman Filter-predicted state trajectory for the state $x_1$. The illustration of Fig. 5(b) is associated with the state $x_2$. In Figs. 5(a) and 5(b) solid lines indicate the true state trajectories, dash-dash and dot-dot lines denote the Carleman linearized estimate and EKF-predicted state trajectories respectively. It can be inferred that the Carleman linearized mean state trajectory shows the better tracking of the true state in contrast to the EKF result.

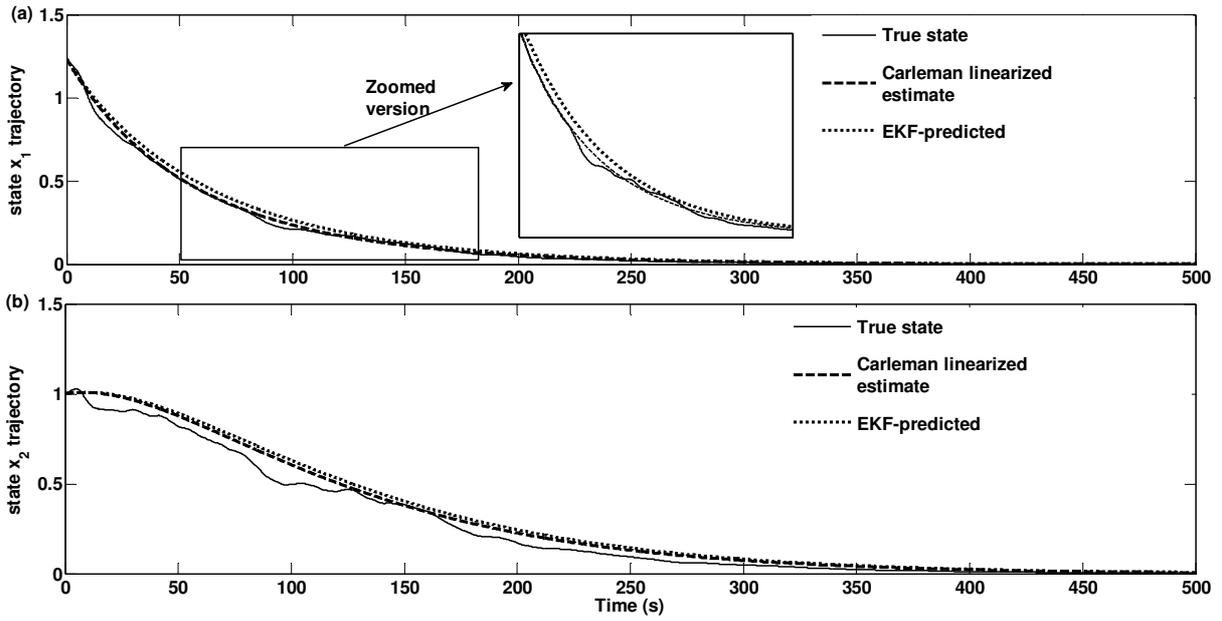

**Fig. 5.** A comparison between true and estimated state trajectories

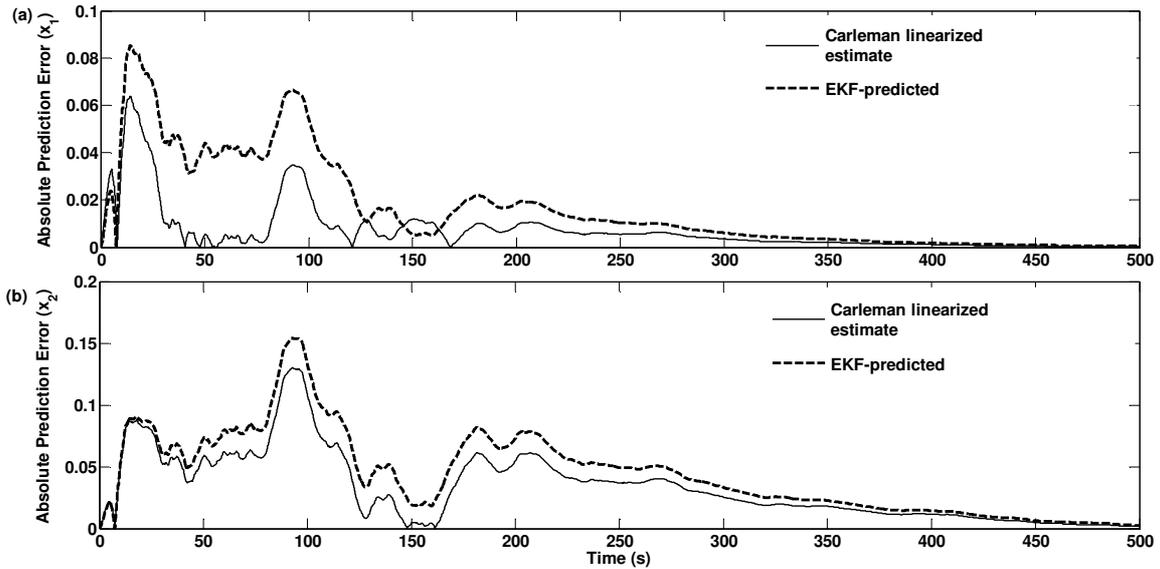

**Fig. 6.** Absolute prediction error comparisons of the both states for the second set of parameters

Figs. 6(a) and 6(b) show a comparison between two absolute prediction error evolutions for the states $x_1$ and $x_2$ respectively. For the given second set of system parameters in Table 3, the maximum absolute prediction error associated with the Carleman linearized estimate is about $0.06$, on the other hand, the absolute maximum error with the EKF-predicted estimate is closer to $0.10$, see Fig. 6(a). Fig. 6(b) illustrates the maximum absolute prediction error associated with the Carleman linearized estimate of the state $x_2$ is less than $0.15$. On the other hand, the maximum absolute error with the EKF-predicted estimate is greater than $0.15$. Fig. 6 reveals the similar inference to that of Fig. 3.

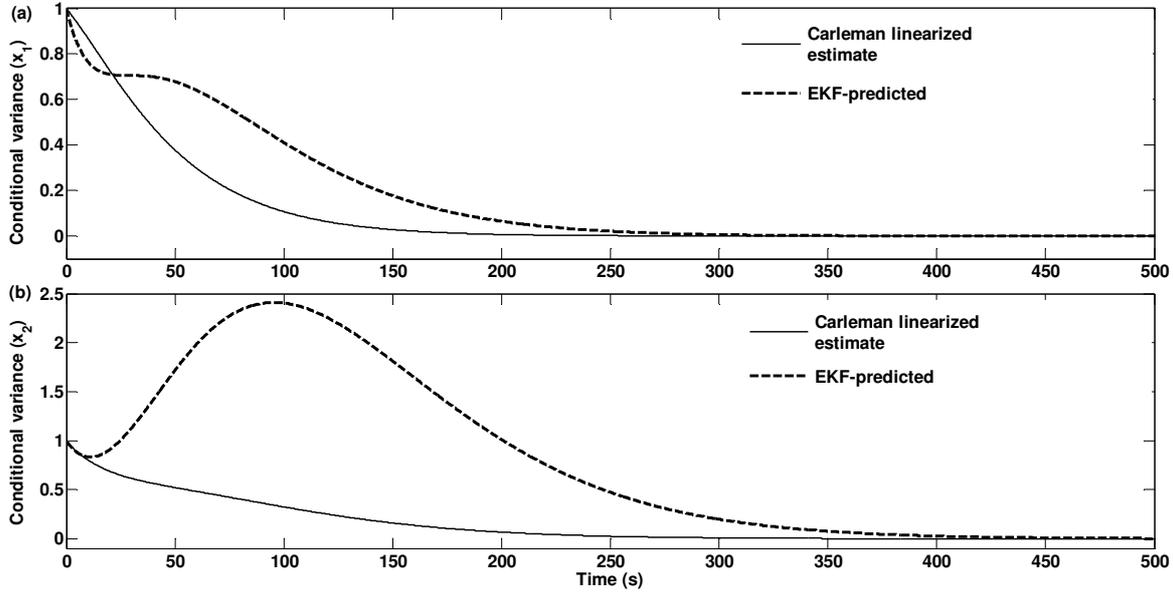

**Fig. 7.** Conditional variance trajectories of the states $x_1$ and $x_2$ for the second set of parameters

Fig. 7 shows conditional variance comparisons resulting from the Carleman linearized estimation and the Extended Kalman Filter-predicted. Fig. 7(a) reveals the less random fluctuation associated with the Carleman linearized estimated state $x_1$ in contrast to the EKF-predicted. This holds for the state $x_2$ as well, see Fig. 7(b).

**Table 4**
Conditional variance of concentrations for the second set of parameters

| Time (s) | Carleman linearized estimated state $P_{x_1}$ | EKF-predicted state $P_{x_1}$ | Carleman linearized estimated state $P_{x_2}$ | EKF predicted state $P_{x_2}$ |
|---|---|---|---|---|
| 0.5 | 0.97 | 0.97 | 0.98 | 0.98 |
| 5 | 0.93 | 0.84 | 0.87 | 0.87 |
| 10 | 0.86 | 0.76 | 0.80 | 0.82 |
| 20 | 0.72 | 0.71 | 0.68 | 0.90 |
| 50 | 0.37 | 0.67 | 0.52 | 1.71 |
| 100 | 0.107 | 0.41 | 0.32 | 2.40 |
| 150 | 0.02 | 0.177 | 0.16 | 1.80 |
| 200 | 0.007 | 0.06 | 0.06 | 1.61 |
| 300 | 0.0005 | 0.0067 | 0.0092 | 0.1992 |
| 400 | 0.00003 | 0.0005 | 0.00109 | 0.0277 |

Table 4 shows the values of conditional variances utilizing the second set of parameters upto 400 seconds. The estimation techniques are the Carleman linearized and the EKF-predicted. The less variance suggests the better estimate. The Carleman linearized estimate is better than that of the EKF-predicted estimate, since the former has the less variance and the latter has the greater variance.

Generally, in chemical reactors, temperature and flow rates are easily available through online measurement systems. However, the reactant and product concentrations are not available directly from the measurement system. They are computed through the on-line analysis. Hence, the need of a better estimation algorithm for the reactant and product concentrations plays a paramount role in the design of system parameters and control strategy. The inaccurate concentration estimation leads to poor design and eventually stability issues arise in chemical reactors. Numerical simulations demonstrated in Figs. 1-7 reveal the following: (i) the Carleman linearized estimate is closer to the true state trajectory (ii) this offers less random fluctuations in the mean trajectory in contrast to the EKF-predicted estimate.

## 5. Conclusion

The main contribution of this paper is to achieve the van de Vusse reactor estimation in the 'Itô framework' via unifying the 'Carleman embedding' and the Fokker-Planck equation. This is the first paper in this direction that adopts a new framework for finding estimation of the van de Vusse reactor.

Secondly, this paper reveals that the Carleman linearized estimate is sharper in contrast to the EKF-predicted estimate. Thus, this paper recommends 'exploring the usefulness of the Carleman linearized estimate as an alternative to the EKF-predicted for non-linear industrial control problems'.

This paper will set new research directions in the sense of 'exploring Carleman embedding for finding sharper and refined estimate and control'. The realization of the Carleman linearized SDE is simpler in contrast to the non-linear as well.

## Appendix

Here, we sketch a proof of the conditional moment evolutions for the Itô non-linear stochastic differential equation with the Carleman linearization order two and the arbitrary size of the state vector.

Making the use of the Fokker-Planck equation and the definition of the differential of the condition expectation of the scalar function of the state vector $\xi_t$, we achieve the exact evolutions of the conditional moment evolution. Consider the double differentiable scalar function $\phi$ of the $n-$dimensional state vector, $\phi: U \to R$, where the phase space $U \subset R^n$. Using the definition of the evolution of conditional expectation of the function $\phi(\xi_t)$, we have

$$d\widehat{\phi}(\xi_t) = \int \phi(\xi) dp(\xi,t|\xi_{t_0},t_0) d\xi = \int \phi(\xi) L(p(\xi,t|\xi_{t_0},t_0)) d\xi dt, \tag{A.1}$$

where $\widehat{\phi}(\xi_t) = E(\phi(\xi_t)|\xi_{t_o},t_0) = \int \phi(\xi) p(\xi,t|\xi_{t_o},t_0) d\xi$ and $E(.)$ is the conditional expectation operator. The term $L(.)$ is the Kolmogorov-Fokker-Planck operator, i.e.

$$L(p(\xi,t|\xi_{t_o},t_0)) = -\sum_i \frac{\partial f_i(\xi,t) p(\xi,t|\xi_{t_0},t_0)}{\partial \xi_i} + \frac{1}{2} \sum_{i,j} \frac{\partial^2 (GG^T)_{i,j}(\xi,t) p(\xi,t|\xi_{t_0},t_0)}{\partial \xi_i \partial \xi_j}. \tag{A.2}$$

After combining (A.1)-(A.2), we have

$$d\widehat{\phi}(\xi_t) = \left\langle \phi, L(p(\xi,t|\xi_{t_0},t_0)) \right\rangle dt. \tag{A.3}$$

Since the Kolmogorov backward operator is an adjoint operator, (A.3) boils down to

$$d\widehat{\phi}(\xi_t) = \left\langle L'\phi, p(\xi,t|\xi_{t_0},t_0) \right\rangle dt, \tag{A.4}$$

where the Kolmogorov backward operator

$$L'(.) = f^T \frac{\partial(.)}{\partial \xi} + \frac{1}{2} tr((G\psi_w G^T)(\xi,t) \frac{\partial^2(.)}{\partial \xi \partial \xi^T}) = \sum_p f_p(t,\xi_t) \frac{\partial(.)}{\partial \xi_p} + \frac{1}{2} \sum_{p,q} (GG^T)_{pq}(t,\xi_t) \frac{\partial^2(.)}{\partial \xi_p \partial \xi_q}. \tag{A.5}$$

As a result of (A.4)-(A.5), we have the following evolution of conditional moment:

$$d\widehat{\phi}(\xi_t) = E\left(\sum_p f_p(\xi_t,t) \frac{\partial \phi(\xi_t)}{\partial \xi_p} + \frac{1}{2}\sum_p (G\psi_w G^T)_{pp}(\xi_t,t) \frac{\partial^2 \phi(\xi_t)}{\partial \xi_p^2} + \sum_{p<q}(G\psi_w G^T)_{pq}(\xi_t,t)\frac{\partial^2 \phi(\xi_t)}{\partial \xi_p \partial \xi_q} \bigg| \xi_{t_0},t_0 \right) dt. \tag{A.6}$$

For the brevity of notations, an alternative set up of (A.6) is

$$d\widehat{\phi}(\xi_t) = \left(\sum_p \widehat{f_p(\xi_t,t) \frac{\partial \phi(\xi_t)}{\partial \xi_p}} + \frac{1}{2}\sum_p \widehat{(G\psi_w G^T)_{pp}(\xi_t,t) \frac{\partial^2 \phi(\xi_t)}{\partial \xi_t^2}} + \sum_{p<q}\widehat{(G\psi_w G^T)_{pq}(\xi_t,t)\frac{\partial^2 \phi(\xi_t)}{\partial \xi_p \partial \xi_q}}\right) dt. \tag{A.7}$$

Thanks to the above (A.7) of the general setting, we weave the specific cases of the conditional moment equation. Thus, the conditional mean and variance evolutions of the Itô stochastic differential equation are

$$d\widehat{\xi}_i = \widehat{f}_i(t,\xi_t) dt, \quad dP_{ij} = (\widehat{\xi_i f_j} - \widehat{\xi}_i \widehat{f}_j + \widehat{f_i \xi_j} - \widehat{f}_i \widehat{\xi}_j + \widehat{(G\varphi_w G^T)_{ij}(\xi_t,t)}) dt. \tag{A.8}$$

Note that $\widehat{\xi}_i = E(\xi_i(t)|\xi_{t_0}, t_0)$. For Itô stochastic differential equation (7), the above-coupled equations, (A.8), boil down to

$$d\widehat{\xi}_t = (A_0(t) + A_t\widehat{\xi}_t)dt, \tag{A.9}$$

$$dP_t = (P_t A_t^T + A_t P_t + G_t G_t^T + G_t \widehat{\xi}_t^T D_t^T + D_t \widehat{\xi}_t G_t^T + D_t \overrightarrow{\widehat{\xi}_t \xi_t^T} D_t^T)dt, \tag{A.10}$$

where $\overrightarrow{\widehat{\xi}_t \xi_t^T} = E(\xi_t \xi_t^T | t_0, \xi_{t_0}) = P_t + \widehat{\xi}_t \widehat{\xi}_t^T$. Equation (A.10) can be rephrased as

$$dP_t = (P_t A_t^T + A_t P_t + G_t G_t^T + G_t \widehat{\xi}_t^T D_t^T + D_t \widehat{\xi}_t G_t^T(t) + D_t P_t D_t^T + D_t \widehat{\xi}_t \widehat{\xi}_t^T D_t^T)dt. \tag{A.11}$$

The last 'five' terms of the right-hand side is the diffusion coefficient contribution. Here, we take a pause and explain the diffusion coefficient contribution to the conditional variance evolution of the bilinear SDE. Note that

$$E(d\xi_t d\xi_t^T | t_0, \xi_{t_0}) = E(((A_0(t) + A_t \xi_t)dt + (G_t + D_t \xi_t)dB_t)((A_0(t)^T + \xi_t^T A_t)dt + dB_t^T(G_t^T + \xi_t^T D_t^T))|t_0, \xi_{t_0})$$

$$= E((G_t + D_t \xi_t)dB_t dB_t^T(G_t^T + \xi_t^T D_t^T)|t_0, \xi_{t_0}). \tag{A.12}$$

After invoking the linearity property of the conditional expectation operator as well as the Itô stochastic differential rule in (A.12), we get

$$E(d\xi_t d\xi_t^T) = (G_t G_t^T + G_t \widehat{\xi}_t^T D_t^T + D_t \widehat{\xi}_t G_t + D_t \overrightarrow{\widehat{\xi}_t \xi_t^T} D_t^T)dt,$$

$$= (G_t G_t^T + G_t \widehat{\xi}_t^T D_t^T + D_t \widehat{\xi}_t G_t^T + D_t P_t D_t^T + D_t \widehat{\xi}_t \widehat{\xi}_t^T D_t^T)dt.$$

Note that the diffusion coefficient associated with the bilinear Itô stochastic differential equation is $G_t G_t^T + G_t \xi_t^T D_t^T + D_t \xi_t G_t^T + D_t \xi_t \xi_t^T D_t^T$, i.e.

$$d\xi_t d\xi_t^T = (G_t G_t^T + G_t \xi_t^T D_t^T + D_t \xi_t G_t^T + D_t \xi_t \xi_t^T D_t^T)dt.$$

Suppose the state vector $\xi_t = (\xi_i) = (x_t^{(k)})$, where $k$ is the Kronecker power, i.e. $1 \leq k \leq N$ $1 \leq i \leq \sum_{1 \leq k \leq N} \binom{n+k-1}{k}$. Note that $\xi_i$ a scalar and the dimensions of the state vectors $x_t$ and $x_t^{(k)}$ are $n$ and $\binom{n+k-1}{k}$ respectively. After excluding redundant state variables from the vector $x_t^{(k)}$, the dimension of the augmented state vector $\xi_t$ reduces to $\sum_{1 \le k \le N} \binom{n+k-1}{k}$ in place of $n^k$.

Introducing the notion of the partitioned vector and matrix format allows to recast them alternatively. Since the Carleman linearization order is two for the stochastic system of the paper, we restrict our discussions to the case, $\xi_t = (x_t^{(k)})$, where $1 \le k \le 2$. Furthermore, the stochastic system is associated with the scalar Brownian motion. Recall (8) and (A.9), (A.11), we have the following conditional mean and variance evolutions in the partitioned vector and matrix format, i.e.

$$d\hat{\xi}_t = d\begin{pmatrix}\hat{x}_t \\ \hat{x}_t^{(2)}\end{pmatrix} = \left(\begin{pmatrix}A_{01}(t) \\ A_{02}(t)\end{pmatrix} + \begin{pmatrix}A_{11}(t) & A_{12}(t) \\ A_{21}(t) & A_{22}(t)\end{pmatrix}\begin{pmatrix}\hat{x}_t \\ \hat{x}_t^{(2)}\end{pmatrix}\right)dt, \quad (A.13)$$

$$dP_{\xi_t} = d\begin{pmatrix}P_{x_t x_t} & P_{x_t x_t^{(2)}} \\ P_{x_t^{(2)} x_t} & P_{x_t^{(2)} x_t^{(2)}}\end{pmatrix} = \left(\begin{pmatrix}P_{x_t x_t} & P_{x_t x_t^{(2)}} \\ P_{x_t^{(2)} x_t} & P_{x_t^{(2)} x_t^{(2)}}\end{pmatrix}\begin{pmatrix}A_{11}^T(t) & A_{21}^T(t) \\ A_{12}^T(t) & A_{22}^T(t)\end{pmatrix} + \begin{pmatrix}A_{11}(t) & A_{12}(t) \\ A_{21}(t) & A_{22}(t)\end{pmatrix}\begin{pmatrix}P_{x_t x_t} & P_{x_t x_t^{(2)}} \\ P_{x_t^{(2)} x_t} & P_{x_t^{(2)} x_t^{(2)}}\end{pmatrix}\right.$$

$$+ \begin{pmatrix}G_1(t)G_1^T(t) & G_1(t)G_2^T(t) \\ G_2(t)G_1^T(t) & G_2(t)G_2^T(t)\end{pmatrix} + \begin{pmatrix}D_{11}(t) & D_{12}(t) \\ D_{21}(t) & D_{22}(t)\end{pmatrix}\begin{pmatrix}\hat{x}_t \\ \hat{x}_t^{(2)}\end{pmatrix}(G_1^T(t) \ G_2^T(t))$$

$$+ \begin{pmatrix}G_1(t) \\ G_2(t)\end{pmatrix}\begin{pmatrix}\hat{x}_t \\ \hat{x}_t^{(2)}\end{pmatrix}^T \begin{pmatrix}D_{11}^T(t) & D_{21}^T(t) \\ D_{12}^T(t) & D_{22}^T(t)\end{pmatrix} + \begin{pmatrix}D_{11}(t) & D_{12}(t) \\ D_{21}(t) & D_{22}(t)\end{pmatrix}\begin{pmatrix}P_{x_t x_t} & P_{x_t x_t^{(2)}} \\ P_{x_t^{(2)} x_t} & P_{x_t^{(2)} x_t^{(2)}}\end{pmatrix}\begin{pmatrix}D_{11}^T(t) & D_{21}^T(t) \\ D_{12}^T(t) & D_{22}^T(t)\end{pmatrix}$$

$$\left. + \begin{pmatrix}D_{11}(t) & D_{12}(t) \\ D_{21}(t) & D_{22}(t)\end{pmatrix}\begin{pmatrix}\hat{x}_t \\ \hat{x}_t^{(2)}\end{pmatrix}(\hat{x}_t^T \ \hat{x}_t^{(2)T})\begin{pmatrix}D_{11}^T(t) & D_{21}^T(t) \\ D_{12}^T(t) & D_{22}^T(t)\end{pmatrix}\right)dt. \quad (A.14)$$

The component-wise description of the above evolutions, (A.13)-(A.14), can be found in (10)-(11d) of the paper. For the brevity of discussions, we omit the repetition.